\documentclass[a4paper,12pt, legno]{article}
\usepackage{amsmath}
\usepackage{amsfonts}
\usepackage{amsthm}
\usepackage{mathrsfs}
\usepackage[T1]{fontenc}
\usepackage{amsmath}
\usepackage{amsfonts}
\usepackage{amsthm}
\usepackage{mathrsfs}
\usepackage{amssymb}
\usepackage{pst-node}
\usepackage{tikz-cd}
\newtheorem{thm}{Theorem}

\newtheorem{cor}{Corollary}

\newtheorem{rem}[]{Remark}
\newcommand{\Rset}{\mathbb{R}}
\newcommand{\Cset}{\mathbb{C}}

\begin{document}
\title{ Which Urbanik class $L_k$, do the hyperbolic and the generalized  logistic characteristic functions belong to?}
\author{Zbigniew J.  Jurek\footnote{Institute of Mathematics, University of Wroc\l aw, Pl. Grunwaldzki 2/4, 50-384 Wroc\l aw, Poland; zjjurek@math.uni.wroc.pl}}
\date{October 15, 2022}
\maketitle

\textbf{Abstract.} Selfdecomposable variables obtained from series of Laplace (double exponential)
variables are objects of this study. We proved that hyperbolic-sine and  hyperbolic-cosine variables are in  the difference of the  Urbanik classes $L_2$ and $L_3$ while generalized logistic variable is  at least in the Urbanik class $L_1$. Hence some ratios of those corresponding  selfdecomposable characteristic functions are  again selfdecomposable.



\medskip
The class of infinitely divisible distributions, ID, plays an important role in the theory of limiting distributions.
 It coincides with limiting distributions of sums of infinitesimal triangular arrays  and is intimately
 connected with  L\'evy stochastic processes. When triangular infinitesimal arrays are obtained from normalized partial
  sums of sequences of independent variables, at a limit, one gets the class, L, of selfdecomposable distributions.
If we have sequences of independent and identically distributed variables we obtain class, $\mathcal{S}$,
of stable distributions  and in particular, Gaussian (normal) distributions.

\noindent For a history of that topic in probability see  Feller (1966), Chapter XVII or  Gnedenko and Kolmogorov (1954),
Sect. 17-19 or Loeve (1963), Sect. 23.

On the other hand, let us also mention  that more recently selfdeomposability appeared in some statistical applications,
in particular, in  models for option pricing in mathematical finance, cf. Carr-Geman-Madan-Yor (2007) or Trabs (2014)),
as well in statistical physics -Ising models, (cf.De Coninck and Jurek (2000) and Jurek (2001)(a)).

\medskip
\textbf{1. Urbanik $L_k, k=0,1,2,...,\infty,$ classes.}

Urbanik (1972) (a summary of results)  and Urbanik (1973) (results with proofs)  introduced and  described  a decreasing
 family of classes,  $L_k, k=0,1, 2,..\infty$, (here $L$ stands for L\'evy's name), of distributions obtained in some schemes of
 limiting procedures, and in   a such way that we have  the following proper inclusions:
\begin{multline}
\mbox{(Gaussian} \subset \mathcal{S}\subset ...\subset L_k \subset L_{k-1}\subset...\subset L_1\subset L_0\equiv L \subset ID,\\
L_{\infty}:=\cap_{k=0}^\infty L_k = \mbox{the smallest closed convolution semigroup} \\ \mbox{containing all stable distributions.}
\end{multline}

Analytically, in terms of \emph{the characteristic function} ( in short: char. f.) , $\phi(t)$, we have the following  characterization
\begin{equation}
\phi(t)\in L_k \ \ \ \mbox{iff}\ \ \forall (0<c<1) \ \phi(t)/\phi(ct)\in L_{k-1},  \ \mbox{where} \ \ L_{-1}:=ID,
\end{equation}
cf. Urbanik (1973), Proposition 1, Theorem 2 and Corollary 1.

In terms of random variables Urbanik classes $L_k$ are described as follows:
\begin{equation}
X\in L_k \  \ \  \ \mbox{iff}\ \ \ \forall (t>0)\ \exists (X_t \in L_{k-1}) \ \  \  X \stackrel{d}{=} e^{-t}X+ X_t, \ \
\end{equation}
(the equality in disribution)  where variables $X_t$ and  $X$ are (stochastically) independent.

\medskip
Recalling the (general) form of the stable  characteristic functions and using the description  (2), we infer that the stable
(in particular, the Gaussian) distribution belong to the  Urbanik  class $L_\infty$.

\medskip
The structural stochastic characterization of variables $X$ (the random integral representation) for the classes $L_k, k=0,1,2,...,$ is the
following
\begin{multline}
X\in L_k  \ \mbox{iff}  \ X=\int_0^\infty e^{-t}dY_X(t), \ Y_X(1)\in L_{k-1} \ \mbox{and} \  \mathbb{E}[\log(1+|Y_X(1)|)]< \infty;   \\
\mbox{or} \\
X\in L_k  \ \mbox{iff}  \ X=\int_0^\infty e^{-t}dZ_X(\frac{t^{k+1}}{k+1}),  \ \mbox{with} \ \mathbb{E}[\log^{k+1}(1+|Z_X(1)|)]< \infty,
\end{multline}
where $(Y_X(t), t\ge0)$ and $(Z_X(t), t\ge 0)$ are L\'evy processes referred to as \emph{the background driving L\'evy process}
(in short: BDLP) that are constructed from the remainders $(X_t, t>0$) in (3). Moreover, to the variable $Y_X(1)$  we refer as \emph{the background driving variable} (in short: BDRV); cf. Jurek and Vervaat (1983) for the class $L_0$; or  Jeanblanc, Yor and Chesney (2009), Proposition 11, p. 597.  For other classes $L_k$, see Jurek (1983)(a). Comp. also Sato (1980).

In the past mostly the class $L_0\equiv L$ of \emph{ selfdecomposable distributions},
(this terminology is justified by the decomposition in (3)), was studied and applied in mathematical finance
or statistical physics.

\medskip
See Jurek and Mason (1993), Chapter 3, and references therein, or Jurek (1983)(b), for the generalization  of Urbanik classes to infinite dimensional  Banach space valued random vectors and the normalization by bounded linear operators.

[For a general conjecture concerning \emph{random integral representations}, see: $www.math.uni.wroc.pl/^{\sim}zjjurek/$]

\medskip
\medskip
\textbf{2. Results and some corollaries.}

In this note we consider primary selfdecomposable distributions that are  obtained as  sums of series of \emph{ double exponential $\eta$} random variables ( also called \emph{Laplace distributions}). We prove
 \begin{thm}
 (a) \ The hyperbolic-sine $\hat{S}$ and hyperbolic-cosine $\hat{C}$ distributions  with the characteristic functions $\phi_{\hat{S}}(t)= \frac{t}{\sinh(t)}$ and $\phi_{\hat{C}}(t)=\frac{1}{\cosh(t)}$, respectively,  belong to the difference $L_2 \setminus L_3$ of Urbanik classes.

 (b) \ The hyperbolic-tangent $\hat{T}$ with characteristic function $\phi_{\hat{T}}(t)=\frac{\tanh(t)}{t}$ and double-exponential $\eta$  distribution with the characteristic function $\phi_{\eta}(t)=1/(1+t^2)$ belong to the difference $L_0\setminus L_1$ of Urbanik classes.

[The same holds for a finite linear combinations of independent hyperbolic-tangent and double exponential variables.]

 (c) \ The logistic distribution $l_{\alpha}$ with characteristic function $\phi_{l_{\alpha}}(t)=|\Gamma(\alpha+it/\pi)/\Gamma(\alpha)|^2$ belongs at least to the Urbanik class $L_1$.
 \end{thm}

 \medskip
Using iteratively the characterization (2) or (3),  the above facts lead to the following corollaries.
\begin{cor}
Since the hyperbolic-sine  $\phi_{\hat{S}}(t)=\frac{t}{\sinh(t)}\in L_2$  therefore:
\begin{multline*}
(i) \noindent  \ \mbox{for any $0<c<1$ functions} \\
 \phi_{\hat{S}_c}(t):= \phi_{\hat{S}}(t)/\phi_{\hat{S}}(ct) = \frac{\sinh(ct)}{c \sinh(t)} \in L_1 \
\mbox{ are  selfdecomposable char. functions};  \qquad \qquad \qquad \qquad \\
(ii) \ \noindent \mbox{for any $0<b, c <1$ functions}  \qquad \qquad \qquad \qquad         \qquad \qquad \qquad \qquad     \qquad \qquad \qquad \qquad        \\
 \ \phi_{\hat{S}_{b,c}}(t):= \phi_{\hat{S}_c}(t)/ \phi_{\hat{S}_c}(bt)
 = \frac{\sinh(ct)\sinh(bt)}{\sinh(t)\sinh(bc t)}\in L_0 \ \mbox { are selfdecomposable ;} \\
(iii) \ \noindent \mbox{for any $0<a,b,c<1$ functions}  \qquad \qquad \qquad \qquad         \qquad \qquad \qquad \qquad     \qquad \qquad \qquad \qquad   \\
\ \ t\to  \phi_{\hat{S}_{b,c}}(t)/\phi_{\hat{S}_{b,c}}(at)
=\frac{\sinh(at)\sinh(bt)\sinh(ct)\sinh(abc t}{\sinh(t)\sinh(abt)\sinh(act)\sinh(bct)}\in ID_{\log}
\\ \ \mbox{are infinitely divisible characteristic functions with finite logarithmic moments.}
\end{multline*}
\end{cor}
The random variable  $\hat{S}_c$,  from (i) in the above corollary,   is called \emph{the Talacko-Zolotarev variable}; cf. You (2022), Sec. 2.2, p.11. It maybe viewed as an innovation variable for the BDRV of hyperbolic-sine as we have:
\begin{cor}
Let  $Y_{\hat{S}}(1)$ be background driving random variable (BDRV) for the hyperbolic-sine  variable $\hat{S}$ and,  for $0<c<1$,  let $\hat{S}_c$ be the Talacko-Zolotarev variable independent of $Y_{\hat{S}}(1)$. Then  we have $$ Y_{\hat{S}}(1)\stackrel{d}{=} cY_{\hat{S}}(1)+ \hat{S}_c, \ \ \mbox{(the equality in distribution.)}$$
 \end{cor}

Similarly as for the case of the hyperbolic-sine, for the hyperbolic-cosine we have the following facts:
\begin{cor}
Since the hyperbolic-cosine  $\phi_{\hat{C}}(t)= \frac{1}{\cosh(t)}\in L_2$ we have that
\begin{multline*}
(i) \ \mbox{for any $0<c<1$,} \\
\ \phi_{\hat{C}_c}(t):=  \phi_{\hat{C}}(t)/\phi_{\hat{C}}(ct)
= \frac{\cosh(ct)}{\cosh(t)} \in L_1 \ \ \mbox{ are  selfdecomposable functions}; \\
\noindent \mbox{For any $0<b<1$ functions}  \qquad \qquad \qquad \qquad         \qquad \qquad \qquad \qquad     \qquad \qquad \qquad \qquad        \\
(ii) \ \phi_{\hat{C}_{b,c}}(t):= \phi_{\hat{C}_c}(t)/ \phi_{\hat{C}_c}(bt)=
\frac{\cosh(ct)\cosh(bt)}{\cosh(t)\cosh(bct)}\in L_0 \ \mbox { are selfdecomposable ;} \\
\noindent \mbox{For any $0<a<1$ functions}  \qquad \qquad \qquad \qquad         \qquad \qquad \qquad \qquad     \qquad \qquad \qquad \qquad   \\
(iii) \ \ t\to   \phi_{\hat{C}_{b,c}}(t)/\phi_{\hat{C}_{b,c}}(at)
=\frac{\cosh(at)\cosh(bt)\cosh(ct)\cosh(abct)}{\cosh(t)\cosh(bct)\cosh(act)\cosh(abt)}\in ID_{\log} \\
\ \mbox{are infinitely dvisible charateristic functions with finite logarithmic moments.}
\end{multline*}
\end{cor}

The same  way as for $\hat{S}_c$ (in Corollary 2)  we may look at the variable $\hat{C}_c$, from above Corollary 3 (i), via  its BDCF
$\psi_{\hat{C}}(t)=\exp[- t\tanh(t)]$.  All in all we get
 \begin{cor}
Let  $Y_{\hat{C}}(1)$ be the background driving random variable (BDRV) for the hyperbolic-cosine  variable $\hat{C}$ and,  for $0<c<1$,  let $\hat{C}_c$ be the variable  in Corollary 3 (i) above,  and independent of $Y_{\hat{C}}(1)$. Then  we have $$ Y_{\hat{C}}(1)\stackrel{d}{=} cY_{\hat{C}}(1)+ \hat{C}_c, \ \ \mbox{(the equality in distribution.)}$$
 \end{cor}

For the logistic characteristic functions we have
\begin{cor}
 (a) As the logistic $l_\alpha \in L_1$, ($\alpha>0$) therefore each $0<c<1$ functions
$$
\Rset \ni t\to |\frac{\Gamma(\alpha +it/\pi)}{\Gamma(\alpha+ ict/\pi)}|^2 \  \ \mbox{are selfdecomposable characteristic functions.}
$$
(b) For $\alpha>0$ and $t\in \Rset$ we have the identity:
$$
\int_0^\infty(\cos(tx)-1)\frac{e^{-\alpha\pi x}}{x \, (1-e^{-\pi x})}dx= \log|\Gamma(\alpha + it/\pi)|-\log \Gamma(\alpha).
$$
\end{cor}

\medskip
\textbf{ 2. Auxiliary facts.}

\medskip
\emph{a). \underline{Selfdecomposable variables among infinitely divisible ones.}}

The classical L\'evy-Khintchine formula gives the description of the infinite divisible random variables $X$ or distributions
 $\mu$ in terms of their characteristic functions. Namely
$$
X\in ID \ \mbox{iff}\ \phi_X(t):=\mathbb{E}[e^{itX}]=\exp[ita-\frac{1}{2}\sigma^2t^2+\int_{\Rset\setminus{\{0\}}}(e^{itx}-1-\frac{itx}{1+x^2})M(dx)],
$$
where the triple $a\in \Rset, \sigma^2\ge0$ and the measure $M$ satisfies the integrability condition:
$\int_{\Rset\setminus{\{0\}}}\min(x^2,1)M(dx)<\infty$,
is uniquely determined. In the sequel, for the simplicity, we will write $X=[a,\sigma^2,M]$, if the above formula holds true.

Recall the following characterization (criterium) for distributions with non-zero L\'evy measures :
\emph{
\begin{multline}
X=[a,\sigma^2,M]\in L_0 \ \mbox{iff} \ M(dx)=k(x)dx, \int_{\Rset\setminus{\{0\}}}k(x)dx=\infty,
\mbox{and function} \\ x \to xk(x) \ \mbox{is non-increasing on both halflines} \ (-\infty, 0)\ \mbox{and} \ (0,\infty);\\
\mbox{and equivalently $(-xk(x))^\prime \ge 0$ is the density of the L\'evy measure} \\ \mbox{of BDRV $Y_X(1)$ in the  first line above (4) }, \ \ \ \ \ \qquad \qquad
\end{multline}
}
cf. Jurek (1997), Corollary 1.1, pp. 95-96 or Jurek and Mason (1993), Theorem 3.4.4, p. 94 or Steutel and van Harn (2004), Theorem 6.12.

\begin{rem}
\emph{ Let $\mathbb{D}^n$ denotes the operator acting on densities of L\'evy measures $k$, (from (5)), of selfdecomposable variables, defined as follows:
$$
(\mathbb{D}^{0}k)(x):= xk(x); \  (\mathbb{D}^1k)(x):= (- x k(x))^\prime ; \ (\mathbb{D}^nk)(x):= \mathbb{D}(\mathbb{D}^{n-1}k)(x),
$$
for $n=2,3,... .$}

\emph{In the examples below, we will be looking for the first $n$ such that $(\mathbb{D}^nk)(x)$ is not a density
of a L\'evy measure of the class $L_0$ distribution, that is, $(\mathbb{D}^nk)(x)$ is an integrable positive function or $(\mathbb{D}^nk)(x)$
assumes negative values.}
\end{rem}

\emph{b). \underline{Series of Laplace (double exponential) variables.}}

Let $\eta$ denotes \emph{Laplace (double exponential) variable} with the probability  density $2^{-1}e^{-|x|}, x\in \Rset$ and $a>0$.
Then $a\eta$  has the characteristic function
\begin{equation}
\phi_{a\eta}(t)=\frac{1}{1+(at)^2}=\exp\int_{\Rset}(\cos(tx)-1)k_{a\eta}(x)dx, \ \ k_{a\eta}(x):= e^{-a^{-1}|x|}/|x|,
\end{equation}
which, for $a=1$, means that  $\eta=[0,0, k_{\eta}] \in ID$ and the corresponding L\'evy measure is equal $M_{\eta}(dx)=k_{\eta}(x)dx$.

\medskip
Since the function $xk_\eta(x)= sign(x)e^{-|x|}$ is non-increasing on both half-lines , by (4), we get that $\eta$ is \emph{selfdecomposable}; in symbols:
$\eta \in L_0$. However, since $(-xk_{\eta}(x))^\prime =e^{-|x|}$ gives finite L\'evy measure thus $\eta \notin L_1$. Consequently, variable
$\eta\in L_0 \setminus L_1.$

\medskip
For independent and identically distributed  Laplace $\eta_k$ variables and  a sequence $\textbf{\underline{a}}:= (a_1,a_2,...)$ of real numbers we have that
\begin{equation}
X(\textbf{\underline{a}}):=\sum_{k=1}^\infty a_k\eta_k< \infty \ \ \mbox{(almost surely)} \ \mbox{iff}\  \ \sum_{k=1}^\infty a_{k}^2 <\infty;
\end{equation}
cf. Jurek (2000) Propositions 1 and 2. Note that without a loss of a generality (as Laplace variables are symmetric)  we may  assume that  $a_k>0$ and sequence $(\underline{a})$  is decreasing  to zero.

Since $L_k$ are closed (in the  weak topology) convolution semigroups (see (2) or (3)) and $\eta\in L_0$ we conclude
that
\begin{equation}
X(\textbf{\underline{a}})=\sum_{k=1}^\infty a_k\eta_k \in L_0, \ \mbox{with} \ \ k_{X(\textbf{\underline{a}})}(x):=\frac{1}{|x|}\sum_{k=1}^\infty e^{-a_k^{-1}|x|},
\end{equation}
and $X(\textbf{\underline{a}})=[0,0, k_{X(\textbf{\underline{a}})}(x)]\in L_0 $.

Hence and  from the  random integral representation (4) for $k=0$, there exist a L\'evy process $(Y_X(t),t\ge 0)$  such that
\begin{multline}
X(\underline{a})=\int_0^\infty e^{-t}dY_{X(\underline{a})}(t), \ \
 h_{{X(\underline{a})}}(x):=\sum_{k=1}^\infty a_k^{-1}e^{-a_k^{-1}|x|}=  (-x  k_{X(\textbf{\underline{a}})}(x))^\prime
\end{multline}
and  $Y_{X(\underline{a})}(1) =[0,0, h_{{X(\underline{a})}}(x)]\in ID_{\log}$ is the background driving random variable (BDRV) for ${X(\underline{a})}$.

\medskip
If $\phi_{X(\underline{a})}(t)$ denotes the characteristic function of $X(\underline{a})$ and $\psi_{X(\underline{a})}(t)$ is the characteristic
of the background driving variable (BDRV) $Y_{X(\underline{a})}(1)$ then
\begin{equation}
\psi_{X(\underline{a})}(t)=\exp[t (\log\phi_{X(\underline{a})}(t))^\prime]= \exp [ t (\phi_{X(\underline{a})}(t))^\prime/\phi_{X(\underline{a})}(t)], t\neq 0;
\end{equation}
cf. Jurek (2001)(b), Corollary 3.  Formulae (9) and (10) give two ways of identifying BDRV $Y_{X(\underline{a})}(1)$.

\begin{rem}
\emph{Series of the form $\sum_{k=1}^\infty a_k^{-1}e^{-a_k^{-1}|x|}$ maybe viewed as a very particular examples of the classical Dirichlet series; cf. Jurek(2000), Section 3 and references therein. These may help to get more explicit examples of variables $X(\underline{a})$. }
\end{rem}

\medskip
\textbf{3. Proofs.}

\medskip
\emph{A).\underline{ Hyperbolic-sine characteristic function $t/\sinh(t)$.}}

From the product representation: \ $\sinh(z)= z \prod_{k=1}^{\infty} (1+\frac{z^2}{k^2\pi^2}), z \in \Cset $,\
taking the sequence $\underline{a}:=  (1/(k\pi)), k=1,2,..$ and putting
$\hat{S}\equiv X(\underline{a})$, by (6),  (7) and (8),  we get the following
\begin{multline*}
\phi_{\hat{S}}(t) = \prod_{k=1}^\infty \frac{1}{1+(t/\pi k)^2}= \frac{t}{\sinh t}, \ k_{\hat{S}}(x)=\frac{1}{|x|}\frac{1}{e^{\pi |x|}-1}
= \frac{e^{-\pi |x|/2}}{2x\sinh(\pi x/2)},
\end{multline*}
that is, $\hat{S}=[0,0,k_{\hat{S}}]\in ID$, with and an infinite L\'evy spectral measure $M_{\hat{S}}(dx):=k_{\hat{S}}(x)dx$.

\medskip
From now on we will employ the procedure described in Remark 1 to the function $k_{\hat{S}}(x)$.

\medskip
\underline{Step 1.} Since for the function $(\mathbb{D}^1k_{\hat{S}})(x)\equiv h_{\hat{S}}(x)$ we have that
$$h_{\hat{S}}(x):= (-xk_{\hat{S}}(x))^\prime  = \pi\frac{e^{\pi|x|}}{(e^{\pi |x|}-1)^2}=\frac{\pi}{4}\frac{1}{\sinh^2(\pi x /2)}=
\frac{\pi}{4}csch^{2}(\pi x/2)>0,
$$
is non-negative, we infer that  the function $x\to x k_{\hat{S}}(x)$ is not increasing on both half-lines which means that \underline{$\hat{S}\in L_0$}.
(Or use the fact from (8)).

Moreover, by (8),  the function $h_{\hat{S}}(x)$ is the density of the L\'evy spectral measure of $Y_{\hat{S}}(1)\in ID_{\log}$,
 where  $(Y_{\hat{S}}(t), t\ge 0)$ is the BDLP for $\hat{S}$.

\medskip
\underline{Step 2.}  Since for  $(\mathbb{D}^2k_{\hat{S}})(x)\equiv g_{\hat{S}}(x)$ we have that the function
\begin{multline}
g_{\hat{S}}(x)= -( x h_{\hat{S}}(x))^\prime = (- \frac{\pi x e^{\pi|x|}}{(e^{\pi|x|}-1)^2})\prime
=  \frac{(- \pi e^{\pi|x|}-\pi^2|x|e^{\pi |x|}}{(e^{\pi|x|}-1)^2} + \frac{2\pi^2|x|e^{2\pi|x|}}{(e^{\pi|x|}-1)^3}\\
=  \frac{\pi e^{\pi|x|}\{e^{\pi|x|}\pi |x|-e^{\pi|x|}+\pi|x|+1\}}{(e^{\pi|x|}-1)^3} = \frac{\pi e^{\pi|x|}}{(e^{\pi|x|}-1)^2} \frac{\pi|x|(e^{\pi|x|}+1)-(e^{\pi|x|}-1)}{(e^{\pi|x|}-1)}\\
= \frac{\pi}{(e^{\pi|x|/2}- e^{-\pi|x|/2})^2}\big[\pi|x|\frac{e^{\pi|x|}+1}{e^{\pi|x|}-1}-1\big]\\= \frac{\pi}{4}csch^{2}(\pi x/2)\,(\pi x\coth(\pi x/2)-1)\ge 0;
\end{multline}
is non-negative, (as the expression in the brackets $\{...\}$ is non-negative; or recall that $x\coth(x)\ge1$), we infer that BDRV  $Y_{\hat{S}}(1) \in L_0$.
Thus,  by  first line in (4), we infer \underline{$\hat{S}\in L_1$}.

Moreover, $g_{\hat{S}}(x)$ it is the  density of L\'evy spectral measure of the background driving  variable $Z_{\hat{S}}(1)\in ID_{\log}$ .

\medskip
\underline{Step 3.}
Again, as before for $(\mathbb{D}^2k_{\hat{S}})(x)\equiv g_{\hat{S}}(x)$, let us  notice that the function
\begin{multline*}
r_{\hat{S}}(x):= - (xg_{\hat{S}}(x))^\prime = - \pi/4 ( x(\pi x \coth(\pi x/2)-1)csch^{2}(\pi x/2))^\prime\\
= \frac{\pi}{8} csch^2(\pi x/2)(2\pi^2 x^2\coth^{2}(\pi x/2) \\ +\pi^2  x^2csch^2(\pi x/2) -6\pi x\coth(\pi x/2)+2)\ge0,
\end{multline*}
is the density of a L\'evy measure of an ID variable.

From the  non-negativity of $r_{\hat{S}}(x)$  we infer that the function
$xg_{\hat{S}}(x)$ is not increasing on both half lines, so  $g_{\hat{S}}$ is a L\'evy function of  $L_0$ variable. Consequently, \underline{$\hat{S}\in L_2$}.

\medskip
\underline{Step 4.}
Finally, putting $(\mathbb{D}^3k_{\hat{S}})(x)\equiv v_{\hat{S}}(x)$  we get that the following function
\begin{multline*}
v_{\hat{S}}(x)= (-\frac{\pi}{8} xr_{\hat{S}}(x))^\prime  \\ = (-  x \frac{\pi}{8}csch^2(\frac{\pi x}{2})(2\pi^2 x^2\coth^{2}(\frac{\pi x}{2})
 +\pi^2  x^2csch^2(\frac{\pi x}{2})-6\pi x\coth(\frac{\pi x}{2})+2)^\prime \\ =
  \frac{\pi}{4} csch^2(\pi x/2)[\pi^3x^3\coth^3(\pi x/2)-6 \pi^2 x^2\coth^2(\pi x/2)-3\pi^2x^2csch^2(\pi x/2) + \\
  \pi x\coth(\pi x/2)(2\pi^2x^2csch^2(\pi x/2)+7)-1].
\end{multline*}
is not positive as , using WolframAlpha, we have that $v_{\hat{S}}(0.9)= -0.0136<0$ !! (or $v_{\hat{S}}(x)<0$ for $0.86<x<1.02$).
Thus it can not be a density function, so  \underline{$\hat{S} \notin L_3$} and $\hat{S}\in L_0\setminus L_3$. This completes a proof of  Theorem  1 (a).

\medskip
\begin{rem}
\emph{ (i) The fact that $\hat{S}\notin L_4$ is also noticed in You (2022) thesis, on p.19, but questions about $L_2$ and $L_3$ were left opened.}

\emph{(ii)  In Talacko (1956) and Zolotarev (1957) one may learn  how these distributions appeared in  statistics and   probability.
Furthermore, all distributions in the above  Corollaries 1  may be viewed as  particular examples of so called Perks' function (ratio of finite sums of exponential functions); cf. Talacko (1956), page 160 or Perks (1932). The same applies to distributions in Corollary 3 below.}

\emph{(iii)
Probability distributions, with the characteristic functions $  \phi_{\hat{S}_c}(t)$ ($0<c<1$) as in  Corollary 2 (i)  are called \emph{Talacko-Zolotarev distributions}.  They are in Urbanik class $L_1$. However, their selfdecomposability (the class $L_0$ property) was already  proved  in  You (2022), Proposition 2.2.1.}
\end{rem}

\medskip
\medskip
\medskip
\medskip
\medskip
\emph{B).\underline{ Hyperbolic-cosine characteristic function $1/\cosh(t)$.} }

Here we proceed along the proof of hyperbolic-sine but we will not use  the mapping  $\mathbb{D}$ from Remark 1
but will keep the same letters for the consecutive densities.

For the hyperbolic-cosine function we have the following  product representation: \   $\cosh(z)=\prod_{k=1}^{\infty} (1+\frac{4z^2}{(2k-1)^2\pi^2}), \ z\in \Cset$.

Taking the sequence $\underline{b}:=  (1/((2k-1)\pi/2)), k=1,2,..$ and  denoting $\hat{C} \equiv X(\underline{b})$ we have
\begin{multline*}
\phi_{\hat{C}}(t) = \prod_{k=1}^\infty \frac{1}{1+(t\pi(2k-1)/2)^2}= \frac{1}{\cosh t} \\  \ \ \
k_{\hat{C}}(x)=\sum_{k=1}^\infty\frac{e^{-\pi/2 (2k-1)|x|}}{|x|}= \frac{e^{-\pi|x|/2}}{|x|(1-e^{-\pi|x|})}= \frac{1}{2|x|\sinh(\pi|x|/2)};
\end{multline*}

\medskip
\underline{Step 1.}
Since the function
\begin{multline*}
h_{\hat{C}}(x)= (-xk_{\hat{C}}(x)^\prime
=\frac{\pi}{2}\frac{e^{\pi|x|/2}+e^{-\pi |x|/2}}{(e^{\pi|x|/2}-e^{-\pi|x|/2})^2}\\=\frac{\pi}{4}\frac{\cosh(\pi|x|/2)}{\sinh^2(\pi|x|/2)}= \frac{\pi}{4}\frac{\cosh(\pi x/2)}{\sinh^2(\pi x/2)}\ge 0,
\end{multline*}
is non-negative we have that \underline{$\hat{C}\in L_0$} and $h_{\hat{C}}(x)$ is the density of the L\'evy measure od the background driving variable $Y_X(1)$.

\medskip
\underline{Step 2.} Since the function
\begin{multline*}
g_{\hat{C}}(x):=-(x h_{\hat{C}}(x))^\prime \\
=\frac{\pi}{8}csch(\frac{\pi x}{2})[\pi x\coth^2(\frac{\pi x}{2})-2\coth(\frac{\pi x}{2})+\pi x csch^2(\frac{\pi x}{2})]\ge 0,
\end{multline*}
is non-negative therefore $Y_X(1)\in L_0$ and thus \underline{$\hat{C}\in L_1$}. Moreover, $g_{\hat{C}}(x)$ is the density of the L\'evy measure for $Y_X(1)$.

\medskip
\underline{Step 3.}
Since  the function
\begin{multline*}
r_{\hat{C}}(x):= - (xg_{\hat{C}}(x))^\prime
=\frac{\pi}{16}csch(\pi x/2)[(\pi x)^2\coth^{3}(\pi x/2)\\ +\coth(\pi x/2)(5(\pi x)^2csch^{2}(\pi x/2)+4)-6\pi x(\coth^2(\pi x/2)-6\pi x csch^2(\pi x/2))]\ge 0
\end{multline*}
is non-negative we have that \underline{$\hat{C}\in L_2$}.

\medskip
\underline{Step 4.}
Finally, since the function
 $v_{\hat{C}}(x):= (-xr_{\hat{C}}(x)^\prime...$ is such that (by WolframAlpha)  $v_{\hat{C}}(2)=-0.346 <0$ we have that \underline{$\hat{C}\notin L_3$}, i.e., that $\hat{C}\in L_0\setminus L_3$, which concludes a proof of Theorem 1 (a).

\medskip
\emph{C). \underline{ Hyperbolic-tangent characteristic function $\tanh(t)/t$.}}

The hyperbolic-tangent $\hat{T}$  has the following L\'evy-Khintchine representation c
$$
\phi_{\hat{T}}(t)=(\tanh t)/t= \exp\int_\Rset(e^{itx}-1-\frac{itx}{1+x^2})\frac{1}{2|x|}[1- \tanh(\pi|x|/4)]dx,
$$
where $k_{\hat{T}}(x):=  \frac{1}{2|x|}[1- \tanh(\pi|x|/4)]$  is the density of L\'evy measures; cf. Jurek and Yor (2004), p.185.

Since the function $xk_{\hat{T}}(x)= 1/2 sign(x)(1-tanh (\pi |x|/4))$ is not increasing on both half lines, by (4), we infer
that \underline{$\hat{T}\in L_0$}.

On the other hand, the function
$
h_{\hat{T}}(x):= (-xk_{\hat{T}}(x))^\prime = \frac{\pi}{8}\frac{1}{\cosh^2(\pi x/4)}
$
is  a density of finite measure L\'evy measure of BDRV $Y_{\hat{T}}(1)$.  Hence again by (4) the hyperbolic tangent \underline{$\hat{T} \notin L_1$}. Thus  $\hat{T}\in L_0 \setminus L_1$, wich proves Theorem 1 (b).

\medskip
\medskip
\emph{D). \underline{Generalized logistic distribution $\beta_\alpha, \alpha>0$.}}

 (a) For the sequence $c_k:= (\pi(\alpha + k-1))^{-1}$ and the variable  $l_\alpha \equiv X(\underline{c})$, by (8),  we have  that the function
 $k_{l_{\alpha}}(x)= \frac{1}{|x|}\frac{e^{-\alpha\pi|x|}}{1-e^{- \pi|x|}} $ is a density of L\'evy measure of $l_\alpha$ variable. Since the function
 \begin{multline*}
   h_{l_\alpha}(x): =(-x  k_{l\alpha}(x))^\prime =
 \frac{\pi e^{-\alpha \pi |x|}(\alpha +(1-\alpha)e^{-\pi|x|}}{(1-e^{-\pi|x|})^2} \\ =
 \frac{\pi}{4}\frac{1}{\sinh^2(\pi|x|/2)} e^{-(\alpha-1)\pi|x|}\{\alpha +(1-\alpha)e^{-\pi|x|}\}\ge 0,
 \end{multline*}
where the non-negativity follows from the fact the expression in $\{...\}$ is non-negative for $\alpha>0$.

By  the criterium (5), $l_{\alpha}\in L_0$ (is selfdecomposable).
[Note that the logistic $l_1$ coincides with hyperbolic-sine function in Section 2 $\textbf{(A)}$.]

Since, (by WolframAlpha) the function $x\to (xh_{l_\alpha})$ is not increasing on both half-lines we infer that $l_{\alpha}$ is in $L_1$.

 \medskip
 Furthermore, using Gradshteyn and Ryzhik (1994), formula \textbf{8.326},1.) we have the L\'evy-Khinchine formula for $l_\alpha$ variable
\begin{multline}
\phi_{l_{\alpha}}(t)= \prod_{k=1}^\infty \frac{1}{1 + (t/(\alpha+k-1)\pi)^2}=|\frac{\Gamma(\alpha + i t/\pi)}{\Gamma(\alpha)}|^2 \\ = \exp \int_{-\infty}^\infty(\cos(tx)-1)k_{l_\alpha}(x) dx.
\end{multline}

\medskip
(b) To have a different approach to the logistic distribution,  let us recall that Euler's beta function $B(x,y)$,  for $ x, y \in \Cset, \Re x>0,\Re y>0$  is defined as
$$
B(x,y):=\int_0^1s^{x-1}(1-s)^{y-1}ds = \int_{-\infty}^\infty \frac{e^{x s}}{(1+e^s)^{x+y}}ds.
$$
For $\alpha>0$,  the random variable $\beta_\alpha$  with the probability density
$$
\frac{1}{B(\alpha, \alpha)} e^{\alpha s}(1+e^s)^{-2\alpha}, \ \mbox{for} \ -\infty<s<\infty.
$$
is  called \emph{a generalized logistic distribution.} Note that
$$
\phi_{\beta_\alpha}(t)= 
\frac{B(\alpha+it,\alpha-it)}{B(\alpha,\alpha)}=|\frac{\Gamma(\alpha+it)}{\Gamma(\alpha)}|^2,
$$
which by (12) means that $\beta_{\alpha}/\pi \stackrel{d}{=}l_\alpha$. So, as before we infer that $\beta_\alpha \in L_1$. Thus Theorem 1 (c) is proved.
Also from (12) we get the part (b) of Corollary 5.

\begin{rem}
\emph{
 Since, $\log \gamma_{\alpha,1}$,  logarithms of gamma variables have characteristic functions $\Gamma(\alpha+it)/\Gamma(\alpha)$,
 cf. Jurek (2021), Example 3.3, we have that $\log \gamma_{\alpha,1}\in L_1$. Thus the above is applicable here as well.}
 \end{rem}

\medskip
\emph{E). \underline{Proofs of Corollaries 2 and 4.}}

\medskip
Recall that for $\phi_X(t)\in L_0$ we define its BDCF (\emph{background driving characteristic function}) as
$$\psi_X(t):= \mathbb{E}[\exp(itY_X(1))]=\exp[t (\log \phi_X(t))^\prime]= \exp[t (\phi_X(t))^\prime / \phi_X(t)], t\neq 0;  $$
where $(Y_X(t), t\ge 0)$ is BDLP;  cf.   (10 ) above and  Jurek (2001), Proposition 3.

Applying the  above for $X=\hat{S}$ and $X=\hat{S}_c$ from Corollary 1 (i), we have
\begin{multline*}
\phi_{\hat{S}}(t)=t/\sinh(t), \ \ \psi_{\hat{S}}(t)= \exp[t (t/\sinh(t))^\prime]=\exp[1- t \coth(t)] \\
\psi_{\hat{S}_c}(t)= \exp[t (\log \phi_{\hat{S}_c}(t))^\prime]=\exp[t (\log(\frac{\sinh(ct)}{c \sinh(t)}))^\prime] \\ =  \exp[ct \coth(ct)-t\coth(t)]=
\exp[(1-t \coth(t))-(1- ct\coth(ct))]\\ = \psi_{\hat{S}}(t)/\psi_{\hat{S}}(ct), \qquad \qquad \qquad \qquad
\end{multline*}
i.e.,  \ $\psi_{\hat{S}}(t) = \psi_{\hat{S}}(ct)\,\psi_{\hat{S_c}}(t)$, or $Y_{\hat{S}}(1)\stackrel{d}{=}c Y_{\hat{S}}(1)+\hat{S}_c$,
which gives Corollary 2.

\medskip
Similarly, for the hyperbolic-cosine  $X=\hat{C}$ and $X=\hat{C}_c$ from Corollary 3 (i), and (10)  we have
\begin{multline*}
\phi_{\hat{C}}(t)=1/\cosh(t); \ \psi_{\hat{C}}(t)=\exp[t(-\log \cosh(t))^\prime]=\exp [- t \tanh(t)], \\ \psi_{\hat{C}_c}(t)= \exp[t (\log \cosh(ct)/\cosh(t))^\prime]
= \exp[t( c\tanh(ct)- \tanh(t))]\\ =\psi_{\hat{C}}(t)/\psi_{\hat{C}}(ct); \ \qquad \qquad \ \mbox{i.e.,} \  \psi_{\hat{C}}(t)= \psi_{\hat{C}}(ct) \,\psi_{\hat{C}_c}(t),\qquad \qquad
\end{multline*}
which completes a proof of Corollary 4.

\medskip
\medskip
\textbf{Acknowledgement.}
This paper was inspired by a short Author's visit to University of Californa, Berkeley, CA,  in March 2022 and discussions
with James Pitman and his Ph D. Student Zhiyi You.  Author is grateful for getting the copy of You's Ph.D. thesis and a copy of
J. Talacko's (1956) paper.

\medskip
\medskip
\textbf{References.}

\medskip
P. Carr, H. Geman, D. B. Madan, M. Yor (2007), \emph{ Self-decomosability and option pricing}, Math. Finance 17 (1), pp.31-57.

\medskip
J. de Coninck and Z. J. Jurek (2000), \emph{Lee-Yang models, selfdecomposability and negative definite functions.}  In: High dimensional probability II; E. Gine, D. M. Mason, J. A. Wellner Editors; \emph{Progress in Prob.} vol. \textbf{47}, Birkh\"auser, pp. 349-367.

\medskip
W. Feller (1966),\emph{An introduction to probability theory and its applications}, John Wiley $\&$ Sons.

\medskip
B.V. Gnedenko and A.N. Kolmogorov (1954), \emph{Limit distributions for sums of independent random variables}, Addison-Wesley, Reading,
Massachusetts.

\medskip
I. S. Gradshteyn and I. M. Ryzhik (1994), \emph{Tables of integrals, series, and products}, Academic Press, San Diego, Fifth Ed.

\medskip
M. Jeanblanc, M. Yor and M. Chesney (2009), \emph{Mathematical methods for financial markets}, Springer 2009.

\medskip
 Z. J. Jurek (1983)(a), The classes $L_{m}(Q)$ of probability measures on Banach spaces. \emph{Bull. Acad. Polon. Sci.},
  \textbf{31}, No 1-2, pp.51-62.

\medskip
Z. J. Jurek (1983)(b),  Limit distributions and one-parameter groups of linear operators on Banach spaces. \emph{J. Multivar. Anal.},
\textbf{13}, no. 4 pp. 578-604

\medskip
Z. J. Jurek (1985), \emph{Relations between the s-selfdecomposable and selfdecomposable measures}, Ann. Probab. 13, 592-608.

\medskip
Z. J. Jurek (1997), \emph{ Selfdecomposability: an exception or a rule?}, Annales Univer. M. Curie-Sk\l odowska, Sectio A, Mathematica, vol.\textbf{51}, pp. 93-107.

\medskip
Z.J. Jurek (2000), \emph{A note on gammma random variables and Dirichlet series}, Stat. Prob. Letters, 49, pp. 387-392.

\medskip
Z. J. Jurek (2001)(a),  \emph{1-D Ising models, geometric random sums and selfdecomposability}, Reports on Math. Physics, 47, pp.21-30.

\medskip
Z. J. Jurek (2001)(b) , \emph{Remarks on the selfdecomposability and new examples},  Demonstration Math. vol. XXXIV (2), pp. 241-250
pp. 85-109.

\medskip
Z. J. Jurek (2021), \emph{On background driving distribution functions (BDDF) for some selfdecopmposable variables},
 Mathematica Applicanda,  vol.49(2), pp. 85-109.

\medskip
Z. J. Jurek and J. D. Mason (1993), \emph{Operator-limit distributions in the probability theory}, J. Wiley$\&$ Sons, New York.

\medskip
Z. J. Jurek and W. Vervaat (1983), \emph{An integral representation for selfdecomposable Banach space valued random variables}, Z. Wahrsch. verw. Gebiete 62, pp. 247-262.

\medskip
Z. J. Jurek  and M. Yor (2004), \emph{Selfdecomposable laws associated with hyperbolic functions}, Probab. Math. Stat., vol. 24, no 1, pp. 41-50.

\medskip
M. Loeve (1963), \emph{Probability theory},Third Edition, D. Van Nostrand Co., Princeton, New Jersey.

\medskip
W. Perks (1932), \emph{On some experiments in the graduation of mortality statistics}, Journal of the Institute of Actuaries, vol. 58, pp. 12-57.

\medskip
K. Sato (1980), \emph{Class L of multivariate distribution and its subclasses},  J. Multivar. Anal. \textbf{10}, pp. 207-232.

\medskip
F. W. Steutel and K. Van Harn (2004), \emph{Infinite divisibility of probability distributions on the real line}, Marcel Dekker, Inc. New York and Basel.

\medskip
J. Talacko (1956), \emph{Perk's distributions and their role in the theory of Wiener's stochastic variables}, In: Trajabos de estadiestica 7.2, pp.159-174.

\medskip
M. Trabs (2014), \emph{Calibration of self-decomposable L\'evy models}, Bernoulli vol. 20, pp. 109-140.

\medskip
K. Urbanik (1972), \emph{ Slowly varying sequences of random variables}, Bull. de L'Acad.  Polon. Sciences, Ser. Math.  Astr. Phys. 20:8 (1972), pp. 679--682.

\medskip
K. Urbanik (1973), \emph{ Limit laws for sequences of normed sums satisfying some stability conditions},
 Proc. 3rd International  Symp. on Multivariate Analysis,   Wright State University, Dayton, OH, USA; June 19-24, 1972,
  Academic  Press, 1973. \newline \noindent[Also on:  www.math.uni.wroc.pl/$^\sim$zjjurek/urb-limitLawsOhio1973.pdf]

\medskip
Z. You (2022), \emph{Some models for dependence in stochastic processes}, Ph. D. Thesis, University of California, Berkeley, Spring 2022.

\medskip
V. M. Zolotarev (1957), \emph{Mellin-Stieltjes transforms in probability theory}, Theor. Probab. Applications, vol. 2, pp. 433-460.
(In Russian).
\end{document}